\documentclass[11pt,epsfig]{article}
\usepackage{amsmath}
\usepackage{amsthm}
 \usepackage{amssymb}
\usepackage{epsfig}
\usepackage{leftidx}
\usepackage{graphicx}
\usepackage{amsfonts}
\usepackage{graphicx, color, epstopdf}
\usepackage{cite}
\usepackage{enumerate}
\usepackage{enumitem}

\newtheorem{theorem}{Theorem}[section]
\newtheorem{lemma}{Lemma}[section]

\newtheorem{corollary}{Corollary}[section]

\newtheorem{remark}{Remark}

\newcommand{\diff}{\triangledown_{\!\tau}}

\newcommand{\defeq}{:=}
\newcommand{\zd}{\,\mathrm{d}}
\newcommand{\abs}[1]{\left|#1\right|}
\newcommand{\absb}[1]{\big|#1\big|}

\newcommand{\bra}[1]{\left(#1\right)}
\newcommand{\brab}[1]{\big(#1\big)}
\newcommand{\braB}[1]{\Big(#1\Big)}
\newcommand{\brat}[1]{(#1)}
\newcommand{\kbra}[1]{\left[#1\right]}
\newcommand{\kbrab}[1]{\big[#1\big]}
\newcommand{\kbraB}[1]{\Big[#1\Big]}

\newcommand{\myinner}[1]{\left\langle#1\right\rangle}
\newcommand{\myinnerb}[1]{\big\langle#1\big\rangle}

\newcommand{\mynorm}[1]{\left\|#1\right\|}
\newcommand{\mynormb}[1]{\big\|#1\big\|}

\newcommand{\mynormt}[1]{\|#1\|}

\def\lan#1{\textcolor{blue}{#1}}

\title{Analysis of adaptive BDF2 scheme for diffusion equations}
\author{Hong-lin Liao\thanks{ORCID 0000-0003-0777-6832; Department of Mathematics,
Nanjing University of Aeronautics and Astronautics,
Nanjing 211106, P. R. China. E-mails: liaohl@nuaa.edu.cn and liaohl@csrc.ac.cn.
He is supported by a grant 1008-56SYAH18037 from NUAA Scientific Research Starting Fund of Introduced Talent.}
\quad Zhimin Zhang\thanks{Beijing Computational Science Research Center, Beijing 100193, P. R. China;
and Department of Mathematics, Wayne State University, Detroit, MI 48202, USA.
E-mails: zmzhang@csrc.ac.cn  and ag7761@wayne.edu.
He is supported partly by the NSFC grant 11871092 and NSAF grant U1530401.}
}
\date{\today}

\begin{document}

\maketitle

\begin{abstract}
The variable two-step backward differentiation formula (BDF2)
is revisited via a new theoretical framework using the positive semi-definiteness of BDF2 convolution kernels
and a class of orthogonal convolution kernels.
We prove that, if the adjacent time-step ratios $r_k:=\tau_k/\tau_{k-1}\le\brat{3+\sqrt{17}}/2\approx3.561$,
the adaptive BDF2 time-stepping scheme for linear reaction-diffusion equations is unconditionally stable and
(maybe, first-order) convergent in the $L^2$ norm. \lan{The second-order temporal convergence can be recovered if
almost all of time-step ratios $r_k\le 1+\sqrt{2}$ or some high-order starting scheme is used.}
Specially, for linear dissipative diffusion problems, the stable BDF2 method preserves
both the energy dissipation law (in the $H^1$ seminorm) and the $L^2$ norm monotonicity at the discrete levels.
An example is included to support our analysis.
\\[1ex]
\emph{Keywords}: linear diffusion equations, adaptive BDF2 scheme,
orthogonal convolution kernels, positive semi-definiteness, stability and convergence\\[1ex]
\emph{AMS subject classifications}: 65M06, 65M12
\end{abstract}

\section{Introduction}
\setcounter{equation}{0}

Adaptive time-stepping strategies are practically useful in capturing the multi-scale behaviors in
many time-dependent differential equations (PDEs). They always require theoretically reliable (stable)
time-stepping methods on arbitrary time meshes, or on general setting of time step-size variations.
For linear and nonlinear parabolic problems, the rigorous numerical analysis of
one-step methods, such us the backward Euler and Crank-Nicolson schemes, may be relatively easy
because they involve only one degree (the current step size $\tau_n$) of freedom.
However, the stability and convergence of multi-step time-stepping approaches with unequal time-steps
would be challenging difficult because they always involve multiple degrees (including the current step $\tau_n$,
the previous step $\tau_{n-1}$ and so on) of freedom.

Due to its strong stability, the variable
two-step backward differentiation formula (BDF2) is practically valuable for stiff or differential-algebraic problems
\cite{Emmrich:2005,GearTu:1974,HairerNorsettWanner:1992,ShampineReichelt:1997,Thomee:2006}.
But the stability and convergence theory remains incomplete so far,
see \cite{Becker:1998,ChenWangYanZhang:2019,Emmrich:2005,LeRoux:1982,Thomee:2006},
even for the simplest linear heat conduction equation $\partial_tu=\Delta u+f$.
In this report, we revisit the BDF2 method from a new point of view
by using the positive semi-definiteness of BDF2 convolution kernels and a novel concept,
namely, discrete orthogonal convolution kernels.
Typically, consider the linear reaction-diffusion problem in a bounded convex domain $\Omega$,
\begin{align}\label{eq: diffusion problem}
\partial_tu-\varepsilon\Delta u&=\kappa(x)u+f(t,x)\quad
\text{for $x\in\Omega$ and $0<t<T$,}
\end{align}
subject to the Dirichlet boundary condition~$u=0$ on the smooth boundary~$\partial\Omega$,
and the initial data $u(0,x)=u_0$ for $x\in\Omega$.
Assume that the diffusive coefficient $\varepsilon>0$ is a constant
and the reaction coefficient $\kappa(x)$ is smooth but bounded by $\kappa^{*}>0$.

Choose the (possibly nonuniform) time levels $0=t_0<t_1<t_2<\cdots<t_N=T$ with the time-step
$\tau_k\defeq t_k-t_{k-1}$ for $1\le k\le N$, and the maximum step size
$\tau\defeq\max_{1\le k\le N}\tau_k$.
For any time sequence $\{v^n\}_{n=0}^N$, denote $\diff v^n:=v^n-v^{n-1}$ and $\partial_{\tau}v^n:=\diff v^n/\tau_n$.
For $k=1,2$, let $\Pi_{n,k}v$ denote the interpolating polynomial of a function~$v$
over $k+1$ nodes $t_{n-k}$, $\cdots$, $t_{n-1}$~and $t_{n}$.
Taking $v^n=v(t_n)$, one can find (for instance, by using
the Lagrange interpolation) that the BDF1 formula $D_1v^n:=\bra{\Pi_{n,1}v}'(t)=\diff  v^n/\tau_n$
for $n\ge1$,
and the BDF2 formula
\begin{align}\label{eq: BDF2 formula}
D_2v^n:=\bra{\Pi_{n,2}v}'(t_n)=&\,\frac{1+2r_{n}}{\tau_n(1+r_{n})}\diff v^{n}
-\frac{r_{n}^2}{\tau_n(1+r_{n})}\diff  v^{n-1}\quad\text{for $n\ge2$,}
\end{align}
where the adjacent time-step ratios
$$r_{k}\defeq\frac{\tau_{k}}{\tau_{k-1}}\;\;\text{for $2\le k\le N$}.$$

Always, one can use the BDF1 scheme, by defining $D_2v^1:=D_1v^1$, to compute the first-level solution $u^1$
because the two-step BDF2 formula needs two starting values and the BDF1 scheme generates a second-order accurate
solution at the first time level. Without losing the generality, we consider only a time-discrete solution, $u^n(x)\approx u(t_n,x)$
for~$x\in\Omega$, is defined by the following adaptive BDF2 time-stepping scheme
\begin{align}\label{eq: time-discrete IBVP}
D_2u^{n}
	=\varepsilon\Delta u^{n}
        +\kappa u^{n}+f^n,\quad\text{for $1\le n\le N$}
\end{align}
where the initial data~$u^0=u_{0}$ and the exterior force $f^n(x)=f(t_n,x)$.
The weak form of the time-discrete problem \eqref{eq: time-discrete IBVP} reads
\begin{align}\label{eq: weak time-discrete IBVP}
\myinnerb{D_2u^{k},w}
	+\varepsilon\myinnerb{\nabla u^{k},\nabla w}
        =\myinnerb{\kappa u^{k},w}+\myinnerb{f^k,w},\quad\text{for $\forall\;w\in  H^1_0(\Omega)$ and $k\ge1$},
\end{align}
where $\myinner{u,w}$ denotes the usual inner product in the space~$L^2(\Omega)$.
Correspondingly,  $\mynorm{\cdot}$ denotes the associated $L^2$ norm
and $\abs{\cdot}_{1}$ is the $H^1$ seminorm. There exists a positive constant $C_{\Omega}$ dependent on the domain $\Omega$
such that $\mynorm{w}\le C_{\Omega}\abs{w}_{1}$ for any $w\in L^2(\Omega)\cap H^1_0(\Omega)$.

Our numerical analysis begins with a new perspective, that is, the BDF2 formula
\eqref{eq: BDF2 formula} is regarded as a discrete convolution summation,
\begin{align}\label{eq: alter BDF2 formula}
D_2v^n=\sum_{k=1}^nb^{(n)}_{n-k}\diff  v^k\quad\text{for $n\ge1$}
\end{align}
where the discrete convolution kernels $b^{(n)}_{n-k}$ are defined by $b^{(1)}_0:=1/\tau_1$, and when $n\ge2$,
\begin{align}\label{eq: BDF2 kernels}
b^{(n)}_0:=\frac{1+2r_n}{\tau_n(1+r_n)},\quad b^{(n)}_1:=-\frac{r_n^2}{\tau_n(1+r_n)}
\quad\text{and}\quad b^{(n)}_j:=0\;\;\text{for $2\leq j\leq n-1$}.
\end{align}
To establish the $L^2$ norm stability and convergence, we introduce a new concept, namely, discrete orthogonal convolution
(DOC) kernels $\big\{\theta_{n-k}^{(n)}\big\}_{k=1}^n$, by a recursive procedure
\begin{align}\label{eq: BDF2 orthogonal procedure}
\theta_{0}^{(n)}:=\frac{1}{b^{(n)}_{0}}\quad\text{and}\quad
\theta_{n-k}^{(n)}:=-\frac{1}{b^{(k)}_{0}}\sum_{j=k+1}^{n}\theta_{n-j}^{(n)}b^{(j)}_{j-k}\quad\text{for $1\leq k\le n-1$.}
\end{align}
Here and hereafter, assume the summation $\sum_{k=i}^{j}$ to be zero and the product $\prod_{k=i}^{j}$ to be one if the index $i>j$.
Obviously, the DOC kernels $\theta_{n-j}^{(n)}$ satisfies the  discrete orthogonal identity
\begin{align}\label{eq: BDF2 orthogonal identity}
\sum_{j=k}^{n}\theta_{n-j}^{(n)}b^{(j)}_{j-k}\equiv\delta_{nk}
\quad\text{for $\forall\;1\leq k\le n$,}
\end{align}
where $\delta_{nk}$ is the Kronecker delta symbol with $\delta_{nk}=0$ if $k\neq n$.
It is to note that the positive semi-definiteness of BDF2 convolution kernels $b^{(n)}_{n-k}$
and the corresponding DOC kernels $\theta_{n-k}^{(n)}$, see Lemmas \ref{lem: BDF2 positive semi-definite}
 and \ref{lem: BDF2orthogonal positive semi-definite}, plays an important role in our numerical analysis.

To make our arguments more clearly, we consider firstly a simple case. Next section focuses on
the linear dissipative parabolic equations \eqref{eq: diffusion problem} with the reaction coefficient $\kappa(x)\le0$.
In the numerical analysis of \eqref{eq: time-discrete IBVP} with $\kappa=0$,
Becker \cite{Becker:1998} proved that, if $0<r_k\le \frac{2+\sqrt{13}}3\approx1.868$,
the discrete solution fulfills,
also see the Thom\'{e}e's classical book \cite[Lemma 10.6]{Thomee:2006},
\begin{align*}
\mynormb{u^{n}}
\le C\exp\bra{C\Gamma_n}\braB{\mynormb{u_{0}}+\sum_{j=1}^{n}\tau_j\mynormb{f^j}}\quad\text{for $n\ge1$,}
\end{align*}
where the quantity $\Gamma_n:=\sum_{k=2}^{n-2}\max\big\{0,r_{k}-r_{k+2}\big\}$ and $C>0$
are dependent on the sequence of step size ratios $r_k$.  To our knowledge,
this type  (but more restrictive) quantity was firstly introduced by Le Roux \cite{LeRoux:1982}
nearly forty years ago; but the quantity $\Gamma_n$ may takes the values of zero, bounded \cite[p.175]{Thomee:2006} or
unbounded \cite[Remark 4.1]{ChenWangYanZhang:2019} at vanishing step sizes by choosing certain step-ratio sequences $\{r_k\}$.
The $L^2$ norm stability estimate by Emmrich \cite{Emmrich:2005}
continues to retain the undesirable prefactor $\exp\bra{C\Gamma_n}$ but improves slightly
the Becker's restriction to $0<r_k\le 1.91$. More recently, with the help of a generalized discrete Gr\"{o}nwall inequality,
Chen \emph{et al} \cite{ChenWangYanZhang:2019}
improved the Becker's estimate by replacing the prefactor $\exp\bra{C\Gamma_n}$ with $\exp\bra{Ct_n}$
but introduced a stronger step-ratio restriction $0<r_k\le 1.53$.
This estimate avoids the worst case of $\Gamma_n$ being unbounded,
but may lose some other approximately ideal situations
with $\Gamma_n=0$.

Note that, the solution of \eqref{eq: diffusion problem} with $\kappa(x)\le0$
satisfies an energy dissipation law
\begin{align}\label{eq: continuous energy law}
\frac{\zd}{\zd t}\bra{\varepsilon \absb{u(t)}_{1}^2+\myinnerb{-\kappa u(t),u(t)}}\le 2\myinnerb{f(t),\partial_tu(t)}\quad \text{for $\kappa\le0$ and $t>0$},
\end{align}
and the following $L^2$ norm estimate
\begin{align}\label{eq: continuous L2 estimate}
\mynormb{u(t)}&\le \mynormb{u_0}+2\int_0^t\mynormb{f(s)}\zd s\quad
\text{for $t\ge0$.}
\end{align}
To the best of our knowledge, the existing $L^2$ norm estimates
for the adaptive BDF2 scheme \eqref{eq: time-discrete IBVP}
always require certain discrete Gr\"{o}nwall inequality so that they are much more
pessimistic than the continuous version \eqref{eq: continuous L2 estimate}, which can be derived
without any continuous Gr\"{o}nwall-Bellman type inequalities.
In section 2, we show in Theorem \ref{thm: enery stability} that the variable BDF2 method has
a discrete energy dissipation law, simulating \eqref{eq: continuous energy law} at the discrete time levels,
when the adjacent time-step ratios $r_k$ satisfy a sufficient condition
\begin{enumerate}[itemindent=1em]
  \item[\textbf{S1}.] $0< r_k\le \brat{3+\sqrt{17}}/2\approx3.561$ for $2\le k\le N$.
\end{enumerate}
Then, with this condition \textbf{S1}, Theorem \ref{thm: L2 norm stability} establishes a novel $L^2$ norm estimate
\begin{align*}
\mynormb{u^{n}}
\le\mynormb{u^{0}}+2\sum_{k=1}^{n}\sum_{j=1}^{k}\theta_{k-j}^{(k)}\mynormb{f^j}
\le\mynormb{u^{0}}+2t_n\max_{1\le j\le n}\mynormb{f^j}\quad\text{for $n\ge1$,}
\end{align*}
which perfectly simulates the continuous estimate \eqref{eq: continuous L2 estimate}.
So the adaptive BDF2 time-stepping method \eqref{eq: time-discrete IBVP} is monotonicity-preserving (in the sense of \cite{HundsdorferRuuthSpiteri:2004}),
unconditionally stable and (maybe, first-order) convergent in the $L^2$ norm.
It is interesting to mention that, the step-ratio condition \textbf{S1}
ensures the A-stability and L-stability of adaptive BDF2 method considering the linear
ODE model equation $y'=\lambda y$ with $\Re\bra{\lambda}\le0$,
see Remark \ref{remark: A stability ODE}.

In Section 3, the stability and convergence of adaptive BDF2 method \eqref{eq: time-discrete IBVP}
is established for the linear diffusion equations \eqref{eq: diffusion problem} with a bounded coefficient $\kappa(x)$.
Theorems \ref{thm: general L2 norm stability} gives
\begin{align*}
\mynormb{u^{n}}\le&\,
2\exp\bra{4\kappa^{*}t_{n-1}}\bigg(\mynormb{u^{0}}+2\sum_{k=1}^{n}\sum_{j=1}^{k}\theta_{k-j}^{(k)}\mynormb{f^j}\bigg)
\quad\text{for $1\le n\le N$,}
\end{align*}
and thus the numerical solution is unconditionally stable with respect to the $L^2$ norm
under the condition \textbf{S1}. As a by-product of the proof for Theorems \ref{thm: general L2 norm stability},
we find that the condition \textbf{S1} also ensures the zero-stability, see Remark \ref{remark: zero stability ODE},
by considering the nonlinear ODE problem $y'=g(t,y)$ with the Lipschitz-continuous  perturbations.
Note that, our stability condition \textbf{S1} updates the classical stability restriction,
$0<r_k<1+\sqrt{2}$, given by Grigorieff \cite{Grigorieff:1983} nearly forty years ago,
also see \cite{CrouzeixLisbona:1984}
and the classical book \cite[Section III.5]{HairerNorsettWanner:1992} by Hairer \emph{et al}.

When the solution varies slowly, we can use the uniform and quasi-uniform  meshes to capture the numerical behavior,
and the adjacent step ratios $r_k$ are always close to 1.
The restriction of adjacent step ratios often takes its effect
in the fast-varying (high gradient) domains, and in the transition regions between the slow-varying and fast-varying domains.
In the ``slow-to-fast" transition regions and fast-varying domains, we need a reduction of time-steps
with the adjacent step ratios $r_k\in(0,1)$,
which is covered by the condition \textbf{S1}.
Actually, the restriction \textbf{S1} limits only the
amplification of time-steps in the ``fast-to-slow" transition regions.
\textbf{S1} says that one can use a series of increasing time-steps with the amplification factor up to 3.561.
Nonetheless, very large time-steps always lead to a loss of numerical precisions
and thus large amplification factors would be rarely
used continuously in practical simulations. So it is reasonable to assume that
the size of $\mathfrak{R}_p$ is very small, where $\mathfrak{R}_p$ is an index set
\begin{align*}
\mathfrak{R}_p:=\Big\{k\,\Big|\,
1+\sqrt{2}\le r_k\le \brat{3+\sqrt{17}}/2\Big.\Big\}.
\end{align*}
Then we prove in Theorem \ref{thm: general L2 norm converegnce} that
the adaptive BDF2 method \eqref{eq: time-discrete IBVP} is second-order convergent in the $L^2$ norm
under the following step-ratio condition
\begin{enumerate}[itemindent=1em]
  \item[\textbf{S2}.] The step ratios $r_k$ are contained in \textbf{S1}, but $\abs{\mathfrak{R}_p}=N_0\ll N$.
\end{enumerate}
This condition seems more theoretically rather than practically. Always, potential users are suggested to
choose all of time-step ratios $r_k\in(0,1+\sqrt{2})$, the Grigorieff's stability restriction, with $N_0=0$ for the second-order accurate computations;
however, the condition \textbf{S2} provides certain redundancy for practical choices of time-steps
in self-adaptive numerical simulations. Numerical tests using random time meshes are presented in Section 4
to support our theoretical results.
We end this article by presenting some concluding remarks in the last section.


\section{Stability analysis for dissipative diffusion problem}
\setcounter{equation}{0}

\subsection{Positive semi-definiteness and energy stability}

We describe firstly a sufficient condition on the adjacent time-step ratios $r_k$
so that the discrete convolution kernels $b_{n-k}^{(n)}$ are positive semi-definite,
which will be essential to the stability and convergence of BDF2 time-stepping scheme.
We consider only certain restriction of each step ratio here.
However, in the adaptive time-stepping process, one can choose the next time-step size $\tau_{m+1}$
(or the step ratio $r_{m+1}$) properly according to the information
from previous time-step ratios $\{r_k\}_{k=2}^{m}$, see more comments in Remark \ref{remark: stability condition}.

\begin{lemma}\label{lem: BDF2 positive semi-definite}
Assume that the adjacent step ratios $r_k$ satisfy \textbf{S1}.
For any real sequence $\{w_k\}_{k=1}^n$ with $n$ entries, it holds that
\begin{align*}
2w_k\sum_{j=1}^kb_{k-j}^{(k)}w_j\ge&\,\frac{r_{k+1}}{1+r_{k+1}}\frac{w_{k}^2}{\tau_k}
-\frac{r_k}{1+r_k}\frac{w_{k-1}^2}{\tau_{k-1}}\quad\text{for $k\ge2$.}
\end{align*}
So the discrete convolution kernels $b_{n-k}^{(n)}$ defined in \eqref{eq: BDF2 kernels}
are positive semi-definite,
\begin{align*}
\sum_{k=1}^nw_k\sum_{j=1}^kb_{k-j}^{(k)}w_j\geq{0}\quad \text{for $n\ge1$.}
\end{align*}
\end{lemma}
\begin{proof}
Applying the inequality $2ab\le a^2+b^2$, one has
\begin{align*}
2w_k\sum_{j=1}^kb_{k-j}^{(k)}w_j=&\,2b^{(k)}_{0}w_k^2+2b^{(k)}_{1}w_{k}w_{k-1}
\ge\brab{2b^{(k)}_{0}+b^{(k)}_{1}}w_{k}^2+b^{(k)}_{1}w_{k-1}^2\\
=&\,\frac{2+4r_k-r_k^2}{\tau_k(1+r_k)}w_{k}^2-\frac{r_k^2}{\tau_k(1+r_k)}w_{k-1}^2
=\frac{2+4r_k-r_k^2}{1+r_k}\frac{w_{k}^2}{\tau_k}
-\frac{r_k}{1+r_k}\frac{w_{k-1}^2}{\tau_{k-1}}\\
=&\,\frac{r_{k+1}}{1+r_{k+1}}\frac{w_{k}^2}{\tau_k}
-\frac{r_k}{1+r_k}\frac{w_{k-1}^2}{\tau_{k-1}}
+\braB{\frac{2+4r_k-r_k^2}{1+r_k}-\frac{r_{k+1}}{1+r_{k+1}}}\frac{w_{k}^2}{\tau_k}
\end{align*}
for $k\ge2$.
If the time-step ratios $0<r_k\le r_s$, where $r_s=\frac{3+\sqrt{17}}{2}$ is the positive root of the equation $2+3r_s-r_s^2=0$,
then the following inequality holds
$$\frac{2+4r_k-r_k^2}{1+r_k}\ge\frac{r_s}{1+r_s}\ge\frac{r_{k+1}}{1+r_{k+1}}\quad\text{for $k\ge2$}.$$
Actually, denote $h(x):=\frac{2+4x-x^2}{1+x}$ and then
$h'(x)=(1+x)^{-2}(x+1+\sqrt{3})\brab{\sqrt{3}-1-x}.$
Consider two cases: (i) If $0<x\le \sqrt{3}-1$, then $h'(x)\ge0$. So $h(r_k)\ge h(0)=2>\frac{r_s}{1+r_s}$.
(ii) If $\sqrt{3}-1<x\le r_s$, then $h'(x)\le 0$. So $h(r_k)\ge h(r_s)=\frac{r_s}{1+r_s}$.
Thus it follows that
\begin{align*}
2w_k\sum_{j=1}^kb_{k-j}^{(k)}w_j\ge&\,\frac{r_{k+1}}{1+r_{k+1}}\frac{w_{k}^2}{\tau_k}
-\frac{r_k}{1+r_k}\frac{w_{k-1}^2}{\tau_{k-1}}\quad\text{for $k\ge2$.}
\end{align*}
Therefore, we have
\begin{align*}
2\sum_{k=1}^nw_k\sum_{j=1}^kb_{k-j}^{(k)}w_j
\ge&\,\frac{2}{\tau_1}w_{1}^2+\frac{r_{n+1}}{1+r_{n+1}}\frac{w_{n}^2}{\tau_n}
-\frac{r_2}{1+r_2}\frac{w_{1}^2}{\tau_{1}}\\
\ge&\,\frac{r_{n+1}}{1+r_{n+1}}\frac{w_{n}^2}{\tau_n}+\frac{2+r_2}{1+r_2}\frac{w_{1}^2}{\tau_{1}}\ge0
\quad\text{for $n\ge1$.}
\end{align*}
It completes the proof.
\end{proof}

\begin{remark}\label{remark: stability condition}
Numerical tests on random time meshes in Section 4 suggest that the step ratio condition \textbf{S1} is not necessary.
While, in mathematical manner, the condition \textbf{S1} is also not necessary since
the positive semi-definiteness of discrete convolution kernels $b_{n-k}^{(n)}$
in \eqref{eq: BDF2 kernels} should be determined by the eigenvalues of the following tridiagonal symmetric matrix
\begin{align*}
\mathbf{B}_2:=\left(
 \begin{array}{ccccc}
    2b_0^{(1)} &b_1^{(2)} && &\vspace{0.15cm}\\
     b_1^{(2)} &2b_0^{(2)} & b_1^{(3)}&&\vspace{0.15cm}\\
     &\ddots  &\ddots&\ddots&\vspace{0.1cm}\\
      &&b_1^{(n-1)}&2b_0^{(n-1)}&b_1^{(n)}\vspace{0.15cm}\\
    &&&b_1^{(n)}&2b_0^{(n)}\\
  \end{array}
\right).
\end{align*}
It is seen that, some weaker (maybe, sufficient and necessary) condition for the positive semi-definiteness
of the matrix $\mathbf{B}_2$ would be a certain combination involving all time-step ratios; however,
it is open to us up to now.
\end{remark}

We now consider the energy ($H^1$ seminorm) stability of BDF2 scheme \eqref{eq: time-discrete IBVP} by
defining a (modified) discrete energy $E^k$,
\begin{align}\label{eq: discrete energy}
E^k:=\frac{r_{k+1}}{1+r_{k+1}}\tau_k\mynormb{\partial_\tau u^k}^2+\varepsilon \absb{u^k}_{1}^2+\myinnerb{-\kappa u^{k},u^k}\quad \text{for $\kappa\le0$ and $k\ge1$},
\end{align}
together with the initial energy $E^0:=\varepsilon \absb{u^0}_{1}^2+\myinnerb{-\kappa u^{0},u^0}$.

\begin{theorem}\label{thm: enery stability}
Under the condition \textbf{S1},
the discrete solution $u^n$ of the BDF2 time-stepping scheme \eqref{eq: time-discrete IBVP} with $\kappa\le0$ satisfies
\begin{align}\label{eq: discrete energy law}
	\partial_{\tau}E^k\le 2\myinnerb{f^k,\partial_{\tau}u^k},\quad\text{for $k\ge1$},
\end{align}
which simulates the energy dissipation law \eqref{eq: continuous energy law} numerically.
So the discrete solution is unconditionally stable in the energy norm,
\begin{align*}
	\sqrt{E^n}\le \sqrt{E^0}
+4\varepsilon^{-\frac12}C_{\Omega}\braB{\mynormb{f^1}+\sum_{k=2}^{n}\mynormb{\diff f^{k}}}\quad\text{for $n\ge1$}.
\end{align*}
\end{theorem}
\begin{proof}
 Taking $w=2\diff u^k$ in the weak form \eqref{eq: weak time-discrete IBVP}
 for $k\ge2$, we have
\begin{align*}
2\myinnerb{D_2u^{k},\diff u^k}
	+2\varepsilon\myinnerb{\nabla u^{k},\diff \nabla u^k}
        +2\myinnerb{-\kappa u^{k},\diff u^k}=2\myinnerb{f^k,\diff u^k},\quad\text{for $k\ge2$}.
\end{align*}
Lemma \ref{lem: BDF2 positive semi-definite} gives
\begin{align*}
2\myinnerb{D_2u^{k},\diff u^k}\ge \frac{r_{k+1}}{1+r_{k+1}}\tau_k\mynormb{\partial_\tau u^k}^2
-\frac{r_k}{1+r_k}\tau_{k-1}\mynormb{\partial_\tau u^{k-1}}^2\quad\text{for $k\ge2$}.
\end{align*}
With the help of the inequality $2a(a-b)\ge a^2-b^2$, it is easy to obtain that
\begin{align}\label{thmProof: enery law}
	\diff E^k\le2\myinnerb{f^k,\diff u^k},\quad\text{for $k\ge2$}.
\end{align}
Also, by taking $w=2\diff u^1$ in \eqref{eq: weak time-discrete IBVP} for the case $k=1$, we get
\begin{align*}
2\tau_1\mynormb{\partial_\tau u^1}^2+\varepsilon \absb{u^1}_{1}^2+\myinnerb{-\kappa u^{1},u^1}
\le \varepsilon \absb{u^0}_{1}^2+\myinnerb{-\kappa u^{0},u^0}+2\myinnerb{f^1,\diff u^1},
\end{align*}
which implies that
\begin{align*}
\diff E^1\le2\myinnerb{f^1,\diff u^1}.
\end{align*}
Combining it with the general case \eqref{thmProof: enery law}, one gets the discrete energy dissipation law \eqref{eq: discrete energy law}.
Summing the inequality \eqref{eq: discrete energy law} from $k=1$ to $n$, we have
\begin{align}\label{thmProof: global enery law}
	E^n\le E^0+2\sum_{k=1}^n\myinnerb{f^k,\diff u^k}\quad\text{for $n\ge1$}.
\end{align}
By applying the technique of time summation by parts (cf. \cite[Lemma 2.6]{LiaoSun:2010}) and the Cauchy-Schwarz inequality, we obtain
\begin{align*}
	\sum_{k=1}^n\myinnerb{f^k,\diff u^k}
=&\,\myinnerb{f^n,u^n}-\sum_{k=2}^{n}\myinnerb{\diff f^{k}, u^{k-1}}-\myinnerb{f^1,u^0}\\
\le&\,\mynormb{u^n}\mynormb{f^n}+
\sum_{k=2}^{n}\mynormb{u^{k-1}}\mynormb{\diff f^{k}}+
\mynormb{u^0}\mynormb{f^1}\\
\le&\,\varepsilon^{-\frac12}C_{\Omega}\braB{\sqrt{E^n}\mynormb{f^n}+
\sum_{k=2}^{n}\sqrt{E^{k-1}}\mynormb{\diff f^{k}}+
\sqrt{E^0}\mynormb{f^1}}\quad\text{for $n\ge1$},
\end{align*}
where the Poincar\'{e} inequality has been used. It follows from \eqref{thmProof: global enery law} that
\begin{align*}
	E^n\le E^0+2\varepsilon^{-\frac12}C_{\Omega}\braB{\sqrt{E^n}\mynormb{f^n}+
\sum_{k=2}^{n}\sqrt{E^{k-1}}\mynormb{\diff f^{k}}+
\sqrt{E^0}\mynormb{f^1}}\quad\text{for $n\ge1$}.
\end{align*}
For any finite $n$, choose $n_0$ ($0\leq n_0\le n$) so that $E^{n_0}=\max_{0\leq j\le n}E^j$.
One can take $n=n_0$ in the above inequality and apply the triangle inequality to obtain
\begin{align*}
	E^{n_0}\le& \sqrt{E^0}\sqrt{E^{n_0}}+2\varepsilon^{-\frac12}C_{\Omega}\sqrt{E^{n_0}}\braB{\mynormb{f^{n_0}}+
\sum_{k=2}^{n_0}\mynormb{\diff f^{k}}+\mynormb{f^1}}\\
\le& \sqrt{E^0}\sqrt{E^{n_0}}+4\varepsilon^{-\frac12}C_{\Omega}\sqrt{E^{n_0}}\braB{\mynormb{f^1}+\sum_{k=2}^{n_0}\mynormb{\diff f^{k}}}
\end{align*}
because $f^{n_0}=f^1+\sum_{k=2}^{n_0}\diff f^k$. Thus it follows that
\begin{align*}
	\sqrt{E^n}\le \sqrt{E^{n_0}}\le& \sqrt{E^0}
+4\varepsilon^{-\frac12}C_{\Omega}\braB{\mynormb{f^1}+\sum_{k=2}^{n_0}\mynormb{\diff f^{k}}}\\
\le& \sqrt{E^0}
+4\varepsilon^{-\frac12}C_{\Omega}\braB{\mynormb{f^1}+\sum_{k=2}^{n}\mynormb{\diff f^{k}}}\quad\text{for $n\ge1$}.
\end{align*}
It yields the claimed estimate and completes the proof.
\end{proof}

If the exterior force $f(x,t)$ is zero-valued, the discrete energy law \eqref{eq: discrete energy law} gives
$$E^k\le E^{k-1}\quad\text{for $k\ge1$},$$
so that the variable-step BDF2 scheme \eqref{eq: time-discrete IBVP}
preserves the energy dissipation law at the discrete levels. This property would be important
in simulating the gradient flow problems, cf. \cite{ChenWangYanZhang:2019,LiaoJiZhang:2019,LiaoWangZhang:2019AC} and the references therein.
However,  the energy estimate in Theorem \ref{thm: enery stability} always leads to a
suboptimal $H^1$ seminorm error estimate with respect to both the temporal accuracy and the dependence of
the diffusive coefficient $\varepsilon$ (especially when $\varepsilon$ is small).

\begin{remark}\label{remark: energy stability}
From the computational view of point, the
discrete energy form $E^k$ in \eqref{eq: discrete energy} suggests that small time-steps
(with small step ratios) are necessary to capture the solution behaviors when
$\mynormt{\partial_tu}$ becomes large, and large time-steps (with some big step ratios) are acceptable to accelerate
the time integration when $\mynormt{\partial_t u}$ is small.
\end{remark}

\subsection{Orthogonal convolution kernels and $L^2$ norm stability}

\begin{lemma}\label{lem: BDF2orthogonal positive semi-definite}
If the BDF2 kernels $b_{n-k}^{(n)}$ in \eqref{eq: BDF2 kernels} are positive semi-definite,
the DOC kernels $\theta_{n-j}^{(n)}$ defined in \eqref{eq: BDF2 orthogonal procedure} are also positive semi-definite.
For any real sequence $\{w_j\}_{j=1}^n$, it holds that
\begin{align*}
\sum_{k=1}^nw_k\sum_{j=1}^k\theta_{k-j}^{(k)}w_j\geq{0}\quad \text{for $n\ge1$.}
\end{align*}
\end{lemma}
\begin{proof}
Given any real sequence $\{w_j\}_{j=1}^n$,
one applies the discrete kernels in \eqref{eq: BDF2 kernels} to define another sequence $\{V_j\}_{k=1}^n$ by
\begin{align*}
V_{j}=-\frac{1}{b_{0}^{(j)}}\sum_{\ell=1}^{j-1}b_{j-\ell}^{(j)}V_{\ell}+\frac{w_j}{b_{0}^{(j)}}
\quad\text{for $j\ge1$,}
\end{align*}
or
\begin{align}\label{lemproof: orthogonal-1}
w_j=\sum_{\ell=1}^jb_{j-\ell}^{(j)}V_{\ell}\quad\text{for $j\ge1$.}
\end{align}
Multiplying both sides of the above equality \eqref{lemproof: orthogonal-1} by the DOC kernels
$\theta_{k-j}^{(k)}$, and summing $j$ from $1$ to $k$,
we find
\begin{align}\label{lemproof: orthogonal-2}
\sum_{j=1}^{k}\theta_{k-j}^{(k)}w_j=&\,\sum_{j=1}^{k}\theta_{k-j}^{(k)}
\sum_{\ell=1}^jb_{j-\ell}^{(j)}V_{\ell}
=\sum_{\ell=1}^kV_{\ell}\sum_{j=\ell}^{k}\theta_{k-j}^{(k)}b_{j-\ell}^{(j)}=V_k\quad\text{for $k\ge1$,}
\end{align}
where the summation order has been exchanged in the second equality
and the orthogonal identity \eqref{eq: BDF2 orthogonal identity} was used in the third one.
Thus it follows from \eqref{lemproof: orthogonal-1}-\eqref{lemproof: orthogonal-2} that
\begin{align*}
\sum_{k=1}^nw_k\sum_{j=1}^k\theta_{k-j}^{(k)}w_j
=\sum_{k=1}^nV_k \sum_{\ell=1}^kb_{k-\ell}^{(k)}V_{\ell}\ge0\quad \text{for $n\ge1$,}
\end{align*}
because the BDF2 kernels $b_{n-k}^{(n)}$ are positive semi-definite.
The proof is completed.
\end{proof}

\begin{corollary}\label{corollary: BDF2orthogonal estimate}
The DOC kernels $\theta_{n-j}^{(n)}$ in \eqref{eq: BDF2 orthogonal procedure} fulfill
\begin{align*}
\sum_{j=1}^{n}\theta_{n-j}^{(n)}\equiv\tau_n\quad\text{such that}\quad \sum_{k=1}^{n}\sum_{j=1}^{k}\theta_{k-j}^{(k)}\equiv t_n\quad\text{for $n\ge1$.}
\end{align*}
\end{corollary}
\begin{proof}This proof is similar to that of Lemma \ref{lem: BDF2orthogonal positive semi-definite}.
Taking $v^n=t_n$ in \eqref{eq: alter BDF2 formula}, one can find that
\begin{align*}
1\equiv\sum_{\ell=1}^jb_{j-\ell}^{(j)}\tau_{\ell}\quad\text{for $j\ge1$.}
\end{align*}
Multiplying both sides of the above equality  by the DOC kernels
$\theta_{n-j}^{(n)}$, and summing $j$ from $1$ to $n$,
we apply the orthogonal identity \eqref{eq: BDF2 orthogonal identity} to find
\begin{align*}
\sum_{j=1}^{n}\theta_{n-j}^{(n)}\equiv&\,\sum_{j=1}^{n}\theta_{n-j}^{(n)}
\sum_{\ell=1}^jb_{j-\ell}^{(j)}\tau_{\ell}
=\sum_{\ell=1}^n\tau_{\ell}\sum_{j=\ell}^{n}\theta_{n-j}^{(n)}b_{j-\ell}^{(j)}=\tau_n\quad\text{for $n\ge1$,}
\end{align*}
as desired. The proof is complete.
\end{proof}

\begin{lemma}\label{lem: BDF2 orthogonal formula}
The DOC kernels $\theta_{n-j}^{(n)}$ in \eqref{eq: BDF2 orthogonal procedure} have an explicit formula
\begin{align*}
\theta_{n-k}^{(n)}=\frac{1}{b^{(k)}_{0}}\prod_{i=k+1}^n\frac{r_{i}^2}{1+2r_{i}}
=\frac{\tau_n}{b^{(k)}_{0}\tau_k}\prod_{i=k+1}^n\frac{r_{i}}{1+2r_{i}}\quad\text{ for $1\leq k\le n$}.
\end{align*}
\end{lemma}
\begin{proof}Denote $\hat{\theta}_{n-k}^{(n)}:=\theta_{n-k}^{(n)}b^{(k)}_{0}$ for $1\leq k\le n$.
For $n=1$, the definition \eqref{eq: BDF2 orthogonal procedure}
yields $\hat{\theta}_{0}^{(1)}=1$. For the index $n\ge2$,
we use the definition \eqref{eq: BDF2 orthogonal procedure} and the BDF2 convolution kernels in \eqref{eq: BDF2 kernels} to find
the following recursive procedure
\begin{align}\label{eq: BDF2 auxiliary recursive}
\hat{\theta}_{0}^{(n)}=1\quad\text{and}\quad
\hat{\theta}_{n-k}^{(n)}=-\frac{b^{(k+1)}_{1}}{b^{(k+1)}_{0}}\hat{\theta}_{n-k-1}^{(n)}
=\frac{r_{k+1}^2}{1+2r_{k+1}}\hat{\theta}_{n-k-1}^{(n)}
\quad\text{for $1\leq k\le n-1$.}
\end{align}
Thus a simple induction yields
\begin{align}\label{eq: BDF2 auxiliary orthogonal formula}
\hat{\theta}_{n-k}^{(n)}=\prod_{i=k+1}^n\frac{r_{i}^2}{1+2r_{i}}>0\quad\text{ for $1\leq k\le n$}.
\end{align}
It yields the claimed formula and completes the proof.
\end{proof}

Now we establish the $L^2$ norm stability of the BDF2 scheme \eqref{eq: time-discrete IBVP} for the case $\kappa\le0$.
\begin{theorem}\label{thm: L2 norm stability}
If the BDF2 kernels $b_{n-k}^{(n)}$ in \eqref{eq: BDF2 kernels} are positive semi-definite
(or the sufficient condition \textbf{S1} holds),
the discrete solution $u^n$ of the adaptive BDF2 scheme \eqref{eq: time-discrete IBVP} with the reaction coefficient $\kappa\le0$
is unconditionally stable in the $L^2$ norm,
\begin{align*}
\mynormb{u^{n}}
\le\mynormb{u^{0}}+2\sum_{k=1}^{n}\sum_{j=1}^{k}\theta_{k-j}^{(k)}\mynormb{f^j}
\le\mynormb{u^{0}}+2t_n\max_{1\le j\le n}\mynormb{f^j}\quad\text{for $n\ge1$.}
\end{align*}
Thus the BDF2 scheme \eqref{eq: time-discrete IBVP} is monotonicity-preserving (taking $f\equiv0$)
according to \cite{HundsdorferRuuthSpiteri:2004}.
\end{theorem}
\begin{proof}
Multiplying both sides of the equation \eqref{eq: time-discrete IBVP} by the DOC kernels
$\theta_{k-n}^{(k)}$, and summing $n$ from $1$ to $k$, we find
\begin{align*}
\sum_{j=1}^{k}\theta_{k-j}^{(k)}D_2u^j=\sum_{j=1}^{k}\theta_{k-j}^{(k)}(\varepsilon\Delta+\kappa) u^j+\sum_{j=1}^{k}\theta_{k-j}^{(k)}f^j
\quad\text{for $k\ge1$}.
\end{align*}
Applying the orthogonal identity \eqref{eq: BDF2 orthogonal identity}, one has
\begin{align*}
\sum_{j=1}^{k}\theta_{k-j}^{(k)}D_2u^j=&\,\sum_{j=1}^{k}\theta_{k-j}^{(k)}
\sum_{\ell=1}^jb_{j-\ell}^{(j)}\diff u^{\ell}
=\sum_{\ell=1}^k\diff u^{\ell}\sum_{j=\ell}^{k}\theta_{k-j}^{(k)}b_{j-\ell}^{(j)}=\diff u^{k}
\quad\text{for $k\ge1$,}
\end{align*}
where the summation order has been exchanged in the second equality.
So we have
\begin{align}\label{eq: ModelODE orthogonal form}
\diff u^{k}=\sum_{j=1}^{k}\theta_{k-j}^{(k)}(\varepsilon\Delta+\kappa) u^j
+\sum_{j=1}^{k}\theta_{k-j}^{(k)}f^j\quad\text{for $k\ge1$.}
\end{align}
Making the inner product of the equation \eqref{eq: ModelODE orthogonal form} with $u^k$,
and summing the resulting equality from $k=1$ to $n$, one has
\begin{align*}
\sum_{k=1}^{n}\myinnerb{u^{k},\diff u^{k}}
=&\,\sum_{k=1}^{n}\sum_{j=1}^{k}\myinnerb{u^{k},\theta_{k-j}^{(k)}(\varepsilon\Delta+\kappa) u^j}+
\sum_{k=1}^{n}\sum_{j=1}^{k}\myinnerb{u^{k},\theta_{k-j}^{(k)}f^j}\\
\le&\,\sum_{k=1}^{n}\sum_{j=1}^{k}\myinnerb{u^{k},\theta_{k-j}^{(k)}f^j}\quad\text{for $n\ge1$,}
\end{align*}
where the following inequality derived by Lemma \ref{lem: BDF2orthogonal positive semi-definite}
has been used,
\begin{align*}
\sum_{k=1}^{n}\sum_{j=1}^{k}\myinnerb{u^{k},\theta_{k-j}^{(k)}(\varepsilon\Delta+\kappa) u^j}
=&\,\sum_{k=1}^{n}\sum_{j=1}^{k}\kbraB{-\varepsilon\myinnerb{\nabla u^{k},\theta_{k-j}^{(k)}\nabla u^j}
+\myinnerb{\kappa u^{k},\theta_{k-j}^{(k)} u^j}}\le0.
\end{align*}
Note that, Lemma \ref{lem: BDF2 orthogonal formula} shows that $\theta_{k-j}^{(k)}>0$.
Then the Cauchy-Schwarz inequality  yields
\begin{align*}
\mynormb{u^{n}}^2+\sum_{k=1}^n\mynormb{\diff u^k}^2\le  \mynormb{u^{0}}^2
+2\sum_{k=1}^{n}\mynormb{u^{k}}\sum_{j=1}^{k}\theta_{k-j}^{(k)}\mynormb{f^j}\quad\text{for $n\ge1$.}
\end{align*}
Taking some integer $n_1$ ($0\le n_1\le n$) such that $\mynormb{u^{n_1}}=\max_{0\le k\le n}\mynormb{u^{k}}$. Taking $n:=n_1$ in the above inequality, one gets
\begin{align*}
\mynormb{u^{n_1}}^2\le  \mynormb{u^{0}}\mynormb{u^{n_1}}+2\mynormb{u^{n_1}}\sum_{k=1}^{n_1}\sum_{j=1}^{k}\theta_{k-j}^{(k)}\mynormb{f^j},
\end{align*}
and thus
\begin{align*}
\mynormb{u^{n}}\le&\,\mynormb{u^{n_1}}\le  \mynormb{u^{0}}
+2\sum_{k=1}^{n_1}\sum_{j=1}^{k}\theta_{k-j}^{(k)}\mynormb{f^j}
\le  \mynormb{u^{0}}+2\sum_{k=1}^{n}\sum_{j=1}^{k}\theta_{k-j}^{(k)}\mynormb{f^j}\quad\text{for $n\ge1$.}
\end{align*}
The claimed second estimate follows from Corollary \ref{corollary: BDF2orthogonal estimate} immediately.
\end{proof}

On the uniform mesh with $r_k\equiv1$, Lemma \ref{lem: BDF2 orthogonal formula} yields
\begin{align*}
\theta_{k-1}^{(k)}=\frac{\tau}{3^{k-1}}\quad\text{and}\quad\theta_{k-j}^{(k)}=\frac{2\tau}{3^{k-j+1}}\quad\text{ for $2\leq j\le k$}.
\end{align*}
So Theorem \ref{thm: L2 norm stability} gives (by exchanging the summation order)
\begin{align*}
\mynormb{u^{n}}\le&\,
\mynormb{u^{0}}+2\mynormb{f^1}\sum_{k=1}^{n}\theta_{k-1}^{(k)}+2\sum_{j=2}^{n}\mynormb{f^j}\sum_{k=j}^{n}\theta_{k-j}^{(k)}\\
=&\,
\mynormb{u^{0}}+3\tau\braB{1-\frac{1}{3^{n}}}\mynormb{f^1}+2\tau\sum_{j=2}^{n}\braB{1-\frac{1}{3^{n-j+1}}}\mynormb{f^j}\quad\text{for $n\ge1$,}
\end{align*}
which recovers our previous result in \cite[Lemma 3.2]{LiaoLyuVong:2018} with a slightly different constant for the term
$\mynorm{f^1}$.
Also, it directly leads to the estimate (1.61) in \cite[Theorem 1.7]{Thomee:2006}. On the other hand,
Theorem \ref{thm: L2 norm stability} recovers the solution estimate \eqref{eq: continuous L2 estimate}
under a mild restriction \textbf{S1}, and essentially improves the existing $L^2$ norm estimates
including the classical one \cite[Lemma 10.6]{Thomee:2006}.
Also, no any discrete Gr\"{o}nwall inequalities have been used
in our $L^2$ norm estimate and
no any restrictions of maximum time-step size are required.

\begin{remark}\label{remark: A stability ODE}
Consider the BDF2 scheme $D_2y^n=\lambda y^n$ for solving the ODE model
$y'=\lambda y$ with $\Re(\lambda)\le0$. Reminding the inequality
$2\Re(\bar{y}^k\diff y^k)\ge |y^k|^2-|y^{k-1}|^2,$
one can follow the proof of Theorem \ref{thm: L2 norm stability} to obtain
$\absb{y^{n}}\le\absb{y^{0}}$ for $n\ge1$.
So the adaptive BDF2 scheme is A-stable under the step ratio condition \textbf{S1}.
Obviously, it is also L-stable considering the limit $\lambda\tau_n\rightarrow-\infty$.
\end{remark}

\lan{\begin{remark}\label{remark: remove BDF1}
Our analysis is fit for any other starting schemes
although the BDF1 scheme is applied here to compute the first-level solution.
To see more clear, multiplying the equation \eqref{eq: time-discrete IBVP} by the DOC kernels
$\theta_{k-n}^{(k)}$ and summing $n$ from $2$ to $k$, we have
\begin{align*}
\sum_{j=2}^{k}\theta_{k-j}^{(k)}D_2u^j=\sum_{j=2}^{k}\theta_{k-j}^{(k)}(\varepsilon\Delta+\kappa) u^j+\sum_{j=2}^{k}\theta_{k-j}^{(k)}f^j
\quad\text{for $k\ge2$}.
\end{align*}
Applying the orthogonal identity \eqref{eq: BDF2 orthogonal identity}, one has
\begin{align*}
\sum_{j=2}^{k}\theta_{k-j}^{(k)}D_2u^j=&\,\sum_{\ell=2}^{k}\diff u^{\ell}\sum_{j=\ell}^k\theta_{k-j}^{(k)}b_{j-\ell}^{(j)}+\theta_{k-2}^{(k)}b_{1}^{(2)}\diff u^{1}
=\diff u^{k}+\theta_{k-2}^{(k)}b_{1}^{(2)}\diff u^{1}
\quad\text{for $k\ge2$,}
\end{align*}
where the summation order has been exchanged in the second equality.
So we have
\begin{align*}
\diff u^{k}=-\theta_{k-2}^{(k)}b_{1}^{(2)}\diff u^{1}+\sum_{j=2}^{k}\theta_{k-j}^{(k)}(\varepsilon\Delta+\kappa) u^j
+\sum_{j=2}^{k}\theta_{k-j}^{(k)}f^j\quad\text{for $k\ge2$.}
\end{align*}
Making the inner product of this equation with $u^k$,
and summing the resulting equality from $k=2$ to $n$,
one applies Lemma \ref{lem: BDF2orthogonal positive semi-definite}
and the Cauchy-Schwarz inequality to find
\begin{align*}
\mynormb{u^{n}}^2-\mynormb{u^{1}}^2
\leq&\,-2\sum_{k=2}^{n}\theta_{k-2}^{(k)}b_{1}^{(2)}\myinnerb{u^{k},\diff u^{1}}
+2\sum_{k=2}^{n}\sum_{j=2}^{k}\myinnerb{u^{k},\theta_{k-j}^{(k)}f^j}\\
\le&\,-2b_{1}^{(2)}\mynormb{\diff u^{1}}\sum_{k=2}^{n}\theta_{k-2}^{(k)}\mynormb{u^{k}}
+2\sum_{k=2}^{n}\mynormb{u^{k}}\sum_{j=2}^{k}\theta_{k-j}^{(k)}\mynormb{f^j}\quad\text{for $n\ge2$.}
\end{align*}
By taking $\mynormb{u^{n_1}}=\max_{1\le k\le n}\mynormb{u^{k}}$, it is easy to get
\begin{align*}
\mynormb{u^{n}}\le&\,\mynormb{u^{n_1}}\le\mynormb{u^{1}}-2b_{1}^{(2)}\mynormb{\diff u^{1}}\sum_{k=2}^{n}\theta_{k-2}^{(k)}
+2\sum_{k=2}^{n}\sum_{j=2}^{k}\theta_{k-j}^{(k)}\mynormb{f^j}\quad\text{for $n\ge2$.}
\end{align*}
From the recursive relationship \eqref{eq: BDF2 auxiliary recursive} with the auxiliary discrete kernels $\hat{\theta}_{k-j}^{(k)}$
in \eqref{eq: BDF2 auxiliary orthogonal formula}, one can obtain that $-b_{1}^{(2)}\theta_{k-2}^{(k)}
=\frac{r_{2}^2}{1+2r_{2}}\hat{\theta}_{k-2}^{(k)}=\hat{\theta}_{k-1}^{(k)}$
and
\begin{align*}
\tau_1\sum_{k=1}^{n}\hat{\theta}_{k-1}^{(k)}=\tau_1+\sum_{k=2}^{n}\tau_1\prod_{i=2}^k\frac{r_{i}^2}{1+2r_{i}}
=\tau_1+\sum_{k=2}^{n}\tau_k\prod_{i=2}^k\frac{r_{i}}{1+2r_{i}}\le t_n
\quad\text{for $1\le n\le N$.}
\end{align*}
Thus by Corollary \ref{corollary: BDF2orthogonal estimate}, we arrive at the following corollary.
\begin{corollary}\label{corollary: stability with any starting scheme}
If the BDF2 kernels $b_{n-k}^{(n)}$ in \eqref{eq: BDF2 kernels} are positive semi-definite
(or the sufficient condition \textbf{S1} holds),
the discrete solution $u^n$ of the adaptive BDF2 scheme \eqref{eq: time-discrete IBVP} with the reaction coefficient $\kappa\le0$
is unconditionally stable in the $L^2$ norm,
\begin{align*}
\mynormb{u^{n}}\le&\,\mynormb{u^{1}}+2\mynormb{\partial_{\tau} u^{1}}\tau_1\sum_{k=2}^{n}\hat{\theta}_{k-1}^{(k)}
+2\sum_{k=2}^{n}\sum_{j=2}^{k}\theta_{k-j}^{(k)}\mynormb{f^j}\\
\le&\,\mynormb{u^{1}}+2t_n\mynormb{\partial_{\tau} u^{1}}
+2t_n\max_{2\le j\le n}\mynormb{f^j}\quad\text{for $n\ge2$.}
\end{align*}
\end{corollary}
Obviously, once the first-level solution $u^1$ and the discrete time derivative $\partial_{\tau} u^{1}$ are second-order accurate,
this estimate yields the second-order accuracy of the adaptive BDF2 scheme
under the step-ratio condition \textbf{S1}, cf. the consistency analysis in subsection 3.2.
\end{remark}}

\section{$L^2$ norm convergence for linear diffusion problems}
\setcounter{equation}{0}

\subsection{Priori estimate}
Now consider the $L^2$ norm priori estimate of the adaptive BDF2 scheme \eqref{eq: time-discrete IBVP}
for a general case $\abs{\kappa(x)}\le \kappa^{*}$.
This situation always needs a discrete Gr\"{o}nwall inequality.

\begin{lemma}\label{lem: discrete Gronwall}
Let $\lambda\ge0$, the time sequences $\{\xi_k\}_{k=0}^N$ and $\{V_k\}_{k=1}^{N}$ be nonnegative. If
$$V_n\le \lambda\sum_{j=1}^{n-1}\tau_jV_j+\sum_{j=0}^{n}\xi_j\quad\text{for $1\leq n\le N$},$$
then it holds that
\begin{align*}
V_n\leq \exp(\lambda t_{n-1})\sum_{j=0}^{n}\xi_j
\quad\text{for\;\; $1\le n\le N$.}
\end{align*}
\end{lemma}
\begin{proof}The proof is standard. Under the induction hypothesis $V_j\leq \exp(\lambda t_{j-1})\sum_{k=0}^{j}\xi_k$ for $1\leq j\leq n-1$,
the desired inequality for the index $n$ follows directly from
$$\lambda\sum_{j=1}^{n-1}\tau_j\exp(\lambda t_{j-1})\leq \lambda\int_0^{t_{n-1}}\exp(\lambda t)\zd t=\exp(\lambda t_{n-1})-1.$$
The principle of induction completes the proof.
\end{proof}
\begin{theorem}\label{thm: general L2 norm stability}
If the BDF2 kernels $b_{n-k}^{(n)}$ in \eqref{eq: BDF2 kernels} are positive semi-definite
(or the sufficient condition \textbf{S1} holds) and the maximum time-step size $\tau\le 1/(4\kappa^{*})$,
the discrete solution $u^n$ of the BDF2 scheme \eqref{eq: time-discrete IBVP}
is unconditionally stable in the $L^2$ norm,
\begin{align}\label{ieq: general L2 norm estimate}
\mynormb{u^{n}}\le&\,
2\exp\bra{4\kappa^{*}t_{n-1}}\bigg(\mynormb{u^{0}}+2\sum_{k=1}^{n}\sum_{j=1}^{k}\theta_{k-j}^{(k)}\mynormb{f^j}\bigg)
\quad\text{for $1\le n\le N$.}
\end{align}
\end{theorem}
\begin{proof}
We can start from \eqref{eq: ModelODE orthogonal form} in the proof of Theorem \ref{thm: L2 norm stability}.
Lemma \ref{lem: BDF2orthogonal positive semi-definite} implies that
$$\sum_{k=1}^{n}\sum_{j=1}^{k}\myinnerb{u^{k},\theta_{k-j}^{(k)}\Delta u^j}\le0.$$
Then, making the inner product of \eqref{eq: ModelODE orthogonal form} by $u^k$,
and summing up the resulting equality from $k=1$ to $n$, one derives that
\begin{align*}
\sum_{k=1}^{n}\myinnerb{u^{k},\diff u^{k}}=&\,\sum_{k=1}^{n}\sum_{j=1}^{k}\myinnerb{u^{k},\theta_{k-j}^{(k)}(\varepsilon\Delta+\kappa) u^j}+
\sum_{k=1}^{n}\sum_{j=1}^{k}\myinnerb{u^{k},\theta_{k-j}^{(k)}f^j}\\
\le&\,\sum_{k=1}^{n}\sum_{j=1}^{k}\myinnerb{ u^{k},\theta_{k-j}^{(k)}\kappa u^j}
+\sum_{k=1}^{n}\sum_{j=1}^{k}\myinnerb{u^{k},\theta_{k-j}^{(k)}f^j}\quad\text{for $1\le n\le N$.}
\end{align*}
Lemma \ref{lem: BDF2 orthogonal formula} implies that the DOC kernels $\theta_{k-j}^{(k)}>0$.
So one applies the Cauchy-Schwarz inequality to find
\begin{align*}
\mynormb{u^{n}}^2\le  \mynormb{u^{0}}^2
+2\kappa^{*}\sum_{k=1}^{n}\mynormb{u^{k}}\sum_{j=1}^{k}\theta_{k-j}^{(k)}\mynormb{u^{j}}
+2\sum_{k=1}^{n}\mynormb{u^{k}}\sum_{j=1}^{k}\theta_{k-j}^{(k)}\mynormb{f^j}\quad\text{for $1\le n\le N$.}
\end{align*}
Choosing some integer $n_2$ ($0\le n_2\le n$) such that
$\mynormb{u^{n_2}}=\max_{0\le k\le n}\mynormb{u^{k}}.$
Then, taking $n:=n_2$ in the above inequality, one gets
\begin{align*}
\mynormb{u^{n_2}}^2\le&\,  \mynormb{u^{0}}\mynormb{u^{n_2}}
+2\kappa^{*}\mynormb{u^{n_2}}\sum_{k=1}^{n_2}\mynormb{u^{k}}\sum_{j=1}^{k}\theta_{k-j}^{(k)}
+2\mynormb{u^{n_2}}\sum_{k=1}^{n_2}\sum_{j=1}^{k}\theta_{k-j}^{(k)}\mynormb{f^j}.
\end{align*}
With the help of Corollary \ref{corollary: BDF2orthogonal estimate}, it follows that
\begin{align*}
\mynormb{u^{n}}\le&\,\mynormb{u^{n_2}}\le  \mynormb{u^{0}}
+2\kappa^{*}\sum_{k=1}^{n_2}\tau_k\mynormb{u^{k}}+2\sum_{k=1}^{n_2}\sum_{j=1}^{k}\theta_{k-j}^{(k)}\mynormb{f^j}\\
\le&\,  \mynormb{u^{0}}
+2\kappa^{*}\sum_{k=1}^{n}\tau_k\mynormb{u^{k}}+2\sum_{k=1}^{n}\sum_{j=1}^{k}\theta_{k-j}^{(k)}\mynormb{f^j}
\quad\text{for $1\le n\le N$.}
\end{align*}
Setting the maximum time-step size $\tau\le 1/(4\kappa^{*})$, one has
\begin{align*}
\mynormb{u^{n}}\le&\,  2\mynormb{u^{0}}
+4\kappa^{*}\sum_{k=1}^{n-1}\tau_k\mynormb{u^{k}}+4\sum_{k=1}^{n}\sum_{j=1}^{k}\theta_{k-j}^{(k)}\mynormb{f^j}
\quad\text{for $1\le n\le N$.}
\end{align*}
Lemma \ref{lem: discrete Gronwall} gives the desired estimate \eqref{ieq: general L2 norm estimate} and completes the proof.
\end{proof}

\begin{remark}\label{remark: zero stability ODE}
Let $g(t,y)$ be a Lipschitz-continuous nonlinear function with the Lipschitz constant $L_g>0$.
Apply the BDF2 time-stepping scheme $D_2y^n=g(t_n,y^n)$ to the nonlinear ODE model
$y'=g(t,y)$ for $0<t\le T$. Assuming the perturbed solution $\bar{y}^n$
solves $D_2\bar{y}^n=g(t_n,\bar{y}^n)+\varepsilon^n$, we can follow the proof of Theorem
\ref{thm: general L2 norm stability} to obtain
\begin{align*}
\absb{y^{n}-\bar{y}^n}
\le&\,2\exp\bra{4L_gt_{n-1}}\braB{\absb{y^{0}-\bar{y}^0}+2t_n\max_{1\le j\le n}\absb{\varepsilon^j}}\quad\text{for $1\le n\le N$.}
\end{align*}
So the BDF2 scheme is zero-stable under the step-ratio scondition \textbf{S1}.
It updates the Grigorieff's stability restriction given in \cite{Grigorieff:1983},
also see the classical book \cite[Section III.5]{HairerNorsettWanner:1992} by Hairer \emph{et al}.
\end{remark}

\subsection{Consistency and convergence}

By using Corollary \ref{corollary: BDF2orthogonal estimate}, the priori estimate \eqref{ieq: general L2 norm estimate} gives
the following estimate
\begin{align*}
\mynormb{u^{n}}\le&\,
2\exp\bra{4\kappa^{*}t_{n-1}}\braB{\mynormb{u^{0}}+2t_n\max_{1\le j\le n}\mynormb{f^j}}
\quad\text{for $1\le n\le N$.}
\end{align*}
This estimate would lead to a loss of time accuracy in error analysis,
cf. Theorem \ref{thm: S1 general L2 norm converegnce} below,
because the BDF1 scheme for the first-level solution $u^1$ is only first-order consistent. Also,
the step ratio condition \textbf{S1} requires $\tau_1/\tau_2\ge\frac{\sqrt{17}-3}{4}\approx0.281$,
and prevents our use of very small initial step size, like $\tau_1=O(\tau^2)$,  to recover the second-order accuracy.
In this sense, the analysis and numerical evidences in the note \cite{Nishikawa:2019}
are inadequate,
although the suggested second-order singly-diagonal implicit Runge-Kutta method
would be reliable to compute the first-level solution.

The loss of accuracy is attributed to the unequal time steps and the associated step ratios
because, as well-known, the uniform BDF2 scheme is globally second-order order accurate.
Actually, in next lemma, the global convolution term combined with the DOC kernels $\theta_{k-j}^{(k)}$ in
the estimate \eqref{ieq: general L2 norm estimate}
is evaluated carefully. To a certain degree, it reveals the error behavior of BDF2 time-stepping with respect to
the unequal time-step sizes.

\begin{lemma}\label{lem: BDF2 global convolution error}
For the consistency error $\eta^j:=D_2u(t_j)-\partial_tu(t_j)$ at $t=t_j$, it holds that
\begin{align*}
\sum_{k=1}^{n}\sum_{j=1}^{k}\theta_{k-j}^{(k)}\mynormb{\eta^j}
\le \tau_1\sum_{k=1}^{n}\hat{\theta}_{k-1}^{(k)}\int_{0}^{t_1}\mynormb{\partial_{tt}u}\zd{t}+
\frac{3}2\sum_{j=1}^{n}\tau_j^2\sum_{k=j}^{n}\hat{\theta}_{k-j}^{(k)}\int_{t_{j-1}}^{t_j}\mynormb{\partial_{ttt}u}\zd{t}
\end{align*}
for $1\le n\le N$, where the discrete kernels $\hat{\theta}_{k-j}^{(k)}$ is given by \eqref{eq: BDF2 auxiliary orthogonal formula}.
\end{lemma}
\begin{proof}
For simplicity, denote
\begin{align*}
G_{t2}^{j}=\int_{t_{j-1}}^{t_j}\mynormb{\partial_{tt}u}\zd{t}\quad\text{and}\quad
G_{t3}^{j}=\int_{t_{j-1}}^{t_j}\mynormb{\partial_{ttt}u}\zd{t}\quad\text{for $j\ge1$.}
\end{align*}
For the case of $j=1$, the consistency error is bounded by
\begin{align*}
\mynormb{\eta^1}\le&\,\frac{1}{\tau_1}\int_{0}^{t_1}\mynormb{\partial_{t}u(t)-\partial_{t}u(t_1)}\zd t
\le\frac{1}{\tau_1}\int_{0}^{t_1}\int_{s}^{t_1}\mynormb{\partial_{tt}u}\zd{t}\zd s\le b_0^{(1)}\tau_1G_{t2}^{1}.
\end{align*}
Then Lemma \ref{lem: BDF2 orthogonal formula} together
with the discrete kernels $\hat{\theta}_{k-j}^{(k)}$ in \eqref{eq: BDF2 auxiliary orthogonal formula} yields
\begin{align}\label{lemProof-globalerror1}
\mynormb{\eta^1}\sum_{k=1}^{n}\theta_{k-1}^{(k)}\le
\tau_1G_{t2}^{1}b_0^{(1)}\sum_{k=1}^{n}\theta_{k-1}^{(k)}\le \tau_1G_{t2}^{1}\sum_{k=1}^{n}\hat{\theta}_{k-1}^{(k)}.
\end{align}
By using the Taylor's expansion formula, one can derive that, also see \cite[Theorem 10.5]{Thomee:2006},
\begin{align*}
\eta^j=&\,-\frac{1+r_j}{2\tau_j}\int_{t_{j-1}}^{t_j}(t-t_{j-1})^2\partial_{ttt}u\zd{t}
+\frac{r_j^2}{2(1+r_j)\tau_j}\int_{t_{j-2}}^{t_j}(t-t_{j-2})^2\partial_{ttt}u\zd{t}\\
=&\,-\frac{1}{2}\brab{b_0^{(j)}-b_1^{(j)}}\int_{t_{j-1}}^{t_j}(t-t_{j-1})^2\partial_{ttt}u\zd{t}
-\frac{1}{2}b_1^{(j)}\int_{t_{j-2}}^{t_{j-1}}(t-t_{j-2})^2\partial_{ttt}u\zd{t}\\
&\,-\frac{1}{2}b_1^{(j)}\int_{t_{j-1}}^{t_j}(t-t_{j-1}+\tau_{j-1})^2\partial_{ttt}u\zd{t}\\
=&\,-\frac{1}{2}b_0^{(j)}\int_{t_{j-1}}^{t_j}(t-t_{j-1})^2\partial_{ttt}u\zd{t}
-\frac{1}{2}b_1^{(j)}\int_{t_{j-2}}^{t_{j-1}}(s-t_{j-2})^2\partial_{ttt}u\zd{t}\\
&\,-\frac{1}{2}b_1^{(j)}\tau_{j-1}\int_{t_{j-1}}^{t_j}\brab{2(t-t_{j-1})+\tau_{j-1}}\partial_{ttt}u\zd{t}\quad\text{for $j\ge2$},
\end{align*}
where the BDF2 convolution kernels \eqref{eq: BDF2 kernels} with $b_0^{(j)}-b_1^{(j)}=(1+r_j)/\tau_j$ have been used. So we
apply  the equality $-b_1^{(j)}/b_0^{(j)}=\frac{r_j^2}{1+2r_j}$ to obtain
\begin{align*}
\mynormb{\eta^j}\le&\,\frac{1}{2}b_0^{(j)}\tau_j^2G_{t3}^{j}
-\frac{1}{2}b_1^{(j)}(2\tau_j+\tau_{j-1})\tau_{j-1}G_{t3}^{j}-\frac{1}{2}b_1^{(j)}\tau_{j-1}^2G_{t3}^{j-1}\\
=&\,\frac{1}{2}b_0^{(j)}\kbra{\tau_j^2G_{t3}^{j}-\frac{b_1^{(j)}}{b_0^{(j)}}(1+2r_j)\tau_{j-1}^2G_{t3}^{j}
-\frac{b_1^{(j)}}{b_0^{(j)}}\tau_{j-1}^2G_{t3}^{j-1}}\\
=&\,b_0^{(j)}\tau_j^2G_{t3}^{j}+\frac{r_j^2\tau_{j-1}^2}{2(1+2r_j)}b_0^{(j)}G_{t3}^{j-1}\quad\text{for $j\ge2$}.
\end{align*}
By using Lemma \ref{lem: BDF2 orthogonal formula} and the recursive formula \eqref{eq: BDF2 auxiliary recursive},
we derive that
\begin{align}\label{lemProof-globalerror2}
\sum_{j=2}^{n}\mynormb{\eta^j}\sum_{k=j}^{n}\theta_{k-j}^{(k)}\le&\,
\sum_{j=2}^{n}\tau_j^2G_{t3}^{j}\sum_{k=j}^{n}\hat{\theta}_{k-j}^{(k)}
+\frac12\sum_{j=2}^{n}\frac{r_j^2\tau_{j-1}^2}{1+2r_j}G_{t3}^{j-1}\sum_{k=j}^{n}\hat{\theta}_{k-j}^{(k)}\nonumber\\
=&\,\sum_{j=2}^{n}\tau_j^2G_{t3}^{j}\sum_{k=j}^{n}\hat{\theta}_{k-j}^{(k)}
+\frac12\sum_{\ell=1}^{n-1}\frac{r_{\ell+1}^2\tau_{\ell}^2G_{t3}^{\ell}}{1+2r_{\ell+1}}\sum_{k=\ell+1}^{n}\hat{\theta}_{k-\ell-1}^{(k)}\nonumber\\
=&\,\sum_{j=2}^{n}\tau_j^2G_{t3}^{j}\sum_{k=j}^{n}\hat{\theta}_{k-j}^{(k)}
+\frac12\sum_{j=1}^{n-1}\tau_{j}^2G_{t3}^{j}\sum_{k=j+1}^{n}\hat{\theta}_{k-j}^{(k)}\nonumber\\
\le&\,\frac32\sum_{j=1}^{n}\tau_j^2G_{t3}^{j}\sum_{k=j}^{n}\hat{\theta}_{k-j}^{(k)}\quad\text{for $2\le n\le N$.}
\end{align}
Then the claimed estimate follows from the following equality
\begin{align*}
\sum_{k=1}^{n}\sum_{j=1}^{k}\theta_{k-j}^{(k)}\mynormb{\eta^j}=\mynormb{\eta^1}\sum_{k=1}^{n}\theta_{k-1}^{(k)}+
\sum_{j=2}^{n}\mynormb{\eta^j}\sum_{k=j}^{n}\theta_{k-j}^{(k)}\,,
\end{align*}
and the above two estimates \eqref{lemProof-globalerror1}-\eqref{lemProof-globalerror2}. It completes the proof.
\end{proof}

To process the error analysis, we apply the formula \eqref{eq: BDF2 auxiliary orthogonal formula} to bound the
terms $\sum_{k=j}^{n}\hat{\theta}_{k-j}^{(k)}$ in
Lemma \ref{lem: BDF2 global convolution error} for $1\le j\le n$.
If all step ratios $r_k$ fulfill the Grigorieff's condition, $0<r_k<1+\sqrt{2}$,
we have $\frac{r_{k}^2}{1+2r_{k}}<1$ for $k\ge2$ and thus
\begin{align*}
\sum_{k=j}^{n}\hat{\theta}_{k-j}^{(k)}=
\sum_{k=j}^{n}\prod_{i=j+1}^k\frac{r_{i}^2}{1+2r_{i}}\le\sum_{k=j}^{n}\braB{\frac{r_{c}^2}{1+2r_{c}}}^{k-j}\le \frac{1+2r_{c}}{1+2r_{c}-r_{c}^2}
\quad\text{for $1\leq j\le n$,}
\end{align*}
where $r_c$ takes the maximum value of all step ratios $r_k$. One has the following extension.
\begin{lemma}\label{lem: BDF2 auxiliary orthogonal estimate}
Consider the discrete kernels
$\hat{\theta}_{k-j}^{(k)}$ in \eqref{eq: BDF2 auxiliary orthogonal formula}. If the step ratios satisfy \textbf{S2}, then
\begin{align}\label{eq: BDF2 auxiliary S2 constant}
\sum_{k=j}^{n}\hat{\theta}_{k-j}^{(k)}\le C_r:=\braB{\frac{\hat{r}_{c}^2}{1+2\hat{r}_{c}}}^{N_0}\frac{1+2r_{c}}{1+2r_{c}-r_{c}^2}
\quad\text{for $1\leq j\le n$,}
\end{align}
where $r_c$ takes the maximum value of all step ratios $r_k\in (0,1+\sqrt{2})$
and $\hat{r}_c$ takes the maximum value of those step ratios $r_k\in\kbrab{1+\sqrt{2},\frac{3+\sqrt{17}}{2}}$ for $2\le k\le N$.
\end{lemma}

Let $\tilde{u}^n:=u(t_n,x)-u^n(x)$ for $n\ge0$. Then the error equation of \eqref{eq: time-discrete IBVP} reads
\begin{align}\label{eq: time-discrete IBVP error}
D_2\tilde{u}^{n}
	=\varepsilon\Delta \tilde{u}^{n}
        +\kappa \tilde{u}^{n}+\eta^n,\quad\text{for $1\le n\le N$}
\end{align}
where the local consistency error $\eta^j=D_2u(t_j)-\partial_tu(t_j)$ for $j\ge1$.
If the step ratios satisfy \textbf{S1} with the maximum time-step size $\tau\le 1/(4\kappa^{*})$,
the priori estimate \eqref{ieq: general L2 norm estimate} in Theorem \ref{thm: general L2 norm stability}  yields
\begin{align}\label{ieq: general L2 norm error estimate 1}
\mynormb{\tilde{u}^{n}}\le&\,
2\exp\bra{4\kappa^{*}t_{n-1}}\braB{\mynormb{\tilde{u}^{0}}+2\sum_{k=1}^{n}\sum_{j=1}^{k}\theta_{k-j}^{(k)}\mynormb{\eta^j}}\quad\text{for $1\le n\le N$.}
\end{align}
\lan{With the help of Lemmas \ref{lem: BDF2 global convolution error}--\ref{lem: BDF2 auxiliary orthogonal estimate},
it is easy to obtain the following result.
\begin{theorem}\label{thm: general L2 norm converegnce}
Let $u(t_n,x)$ and $u^n(x)$ be the solutions of the diffusion problem \eqref{eq: diffusion problem}
and the BDF2 scheme \eqref{eq: time-discrete IBVP}, respectively.
If the step ratio condition \textbf{S2} holds with the maximum time-step size $\tau\le 1/(4\kappa^{*})$,
then the time-discrete solution $u^n$ is convergent in the $L^2$ norm,
\begin{align*}
\mynormb{u(t_n)-u^{n}}\le&\,
2C_r\exp\bra{4\kappa^{*}t_{n-1}}\bigg(\mynormb{u_0-u^{0}}+2\tau_1\!\!\int_{0}^{t_1}\!\!\mynormb{\partial_{tt}u}\zd{t}+
3\sum_{j=1}^{n}\tau_j^2\!\!\int_{t_{j-1}}^{t_j}\mynormb{\partial_{ttt}u}\zd{t}\bigg)
\end{align*}
for $1\le n\le N$, where the mesh-dependent constant $C_r=C_r(N_0,r_c,\hat{r}_c)$ is defined in \eqref{eq: BDF2 auxiliary S2 constant}.
\end{theorem}}

Although large step ratios are allowed in the condition \textbf{S2}, the users are suggested to choose
the Grigorieff's step-ratio restriction $r_k\in(0,1+\sqrt{2})$. In such case, $N_0=0$ and $C_r=\frac{1+2r_{c}}{1+2r_{c}-r_{c}^2}$.
Generally, when the time-step ratios $r_k$ are chosen so that the BDF2 kernels $b_{n-k}^{(n)}$ are positive semi-definite (the condition \textbf{S1} is sufficient),
the series $\sum_{k=1}^{n}\hat{\theta}_{k-1}^{(k)}$ in \eqref{ieq: general L2 norm error estimate 1}
would be unbounded as the step sizes vanish. On the other hand, the solution remains the first-order convergence because
$\tau_1\sum_{k=1}^{n}\hat{\theta}_{k-1}^{(k)}\le t_n$ and, see Corollary \ref{corollary: BDF2orthogonal estimate},
$$\sum_{j=1}^{n}\sum_{k=j}^{n}\tau_j\hat{\theta}_{k-j}^{(k)}
=\sum_{j=1}^{n}\sum_{k=j}^{n}b_0^{(j)}\tau_j\theta_{k-j}^{(k)}\le 2t_n\quad\text{ for $1\le n\le N$}.$$
Then the error estimate \eqref{ieq: general L2 norm error estimate 1} gives the following theorem.
\lan{\begin{theorem}\label{thm: S1 general L2 norm converegnce}
If the BDF2 kernels $b_{n-k}^{(n)}$ in \eqref{eq: BDF2 kernels} are positive semi-definite
(or the sufficient condition \textbf{S1} holds)
and the maximum time-step size $\tau\le 1/(4\kappa^{*})$,
then the solution $u^n$ of BDF2 scheme \eqref{eq: time-discrete IBVP}
 is convergent in the $L^2$ norm in the sense that
\begin{align*}
\mynormb{u(t_n)-u^{n}}\le
2\exp\bra{4\kappa^{*}t_{n-1}}\bigg(\mynormb{u_0-u^{0}}+2t_n\int_{0}^{t_1}\!\!\mynormb{\partial_{tt}u}\zd{t}+
3t_n\max_{1\le j\le n}\tau_j\!\!\int_{t_{j-1}}^{t_j}\!\!\mynormb{\partial_{ttt}u}\zd{t}\bigg),
\end{align*}
for $1\le n\le N$. If the BDF1 scheme in \eqref{eq: time-discrete IBVP} is replaced by some high-order starting scheme,
one can follow the proof of Corollary \ref{corollary: stability with any starting scheme} to derive that
\begin{align*}
\mynormb{u(t_n)-u^{n}}\le2\exp\bra{4\kappa^{*}t_{n-1}}\bigg(&\,\mynormb{u(t_1)-u^{1}}
+2t_n\mynormb{\partial_{\tau}\brab{u(t_1)- u^{1}}}\\
&\,+3t_n\max_{1\le j\le n}\tau_j\int_{t_{j-1}}^{t_j}\mynormb{\partial_{ttt}u}\zd{t}\bigg)\quad\text{for $2\le n\le N$.}
\end{align*}
\end{theorem}
}


\section{Numerical example}
\setcounter{equation}{0}

The nonuniform BDF2 method \eqref{eq: time-discrete IBVP}
together with the Fourier pseudo-spectral in space is applied to solve the heat equation
$\partial_tu=\varepsilon\Delta u+f$ on the space-time domain $(0,2)^2\times(0,1]$.
The exterior force $f$ is chosen so that the equation admits an exact solution $u=e^{-t}\sin 2\pi x\cos2\pi y$.

\begin{table}[!ht]
\begin{center}
\tabcolsep 0pt {Table 4.1 \quad Numerical accuracy on random time mesh for $\varepsilon=1.0$} \vspace*{0.5pt}
\def\temptablewidth{1.0\textwidth}
{\rule{\temptablewidth}{1pt}}
\begin{tabular*}{\temptablewidth}{@{\extracolsep{\fill}}cccccc}
 $N$ &$e(N)$& $\tau$ & Order & $\max r_k$ & $N_1$  \\
\hline
64  &1.56e-02 &1.12e-01&--  &13.74 &3\\
128 &3.24e-03 &6.09e-02&2.27&15.26 &8\\
256 &8.66e-04 &3.27e-02&1.90&32.15 &13\\
512 &1.67e-04 &1.59e-02&2.38&395.6 &26\\
1024&4.45e-05 &7.26e-03&1.91&60.13 &40\\
       \end{tabular*}
       {\rule{\temptablewidth}{1pt}}
       \end{center}\label{table:largeEpsilon}
       \end{table}

 \begin{table}[!ht]
\begin{center}
\tabcolsep 0pt {Table 4.2 \quad Numerical accuracy on random time mesh for $\varepsilon=0.1$} \vspace*{0.5pt}
\def\temptablewidth{1.0\textwidth}
{\rule{\temptablewidth}{1pt}}
\begin{tabular*}{\temptablewidth}{@{\extracolsep{\fill}}cccccc}
 $N$ &$e(N)$& $\tau$ & Order & $\max r_k$ & $N_1$  \\
\hline
64  &9.79e-02 &1.32e-01&--  &10.94 &2\\
128 &2.13e-02 &6.36e-02&2.20&13.62 &7\\
256 &6.06e-03 &3.02e-02&1.82&81.12 &10\\
512 &1.40e-03 &1.49e-02&2.11&604.0 &24\\
1024&3.61e-04 &7.59e-03&1.96&448.2 &53\\
       \end{tabular*}
       {\rule{\temptablewidth}{1pt}}
       \end{center}\label{table:smallEpsilon}
       \end{table}

We consider the arbitrary mesh
with random time-steps
$\tau_{k}=T\epsilon_{k}/S$ for $1\leq k\leq N$,
where $S=\sum_{k=1}^{N}\epsilon_{k}$ and $\epsilon_{k}\in(0,1)$ are random numbers subject to the uniform distribution.
No any special treatments have been used to adjust the time-steps so that some large step ratios appear in our experiments,
see the fifth column in Tables 4.1-4.2.
In each run, the $L^2$ norm error $e(N):=\|u(T)-u^{N}\|$ at the final time $T=1$ is
recorded in Tables 4.1--4.2, in which we also list the maximum time-step size $\tau$, the maximum step ratio
and the number (denote $N_1$ in tables) of time levels with the step ratio $r_k\ge (3+\sqrt{17})/2$.
The experimental rate of convergence
is estimated by $\text{Order}\approx\log_2\bra{e(N)/e(2N)}$. From the current data and more tests not listed here,
we see that the adaptive BDF2 time-stepping is robustly stable and second-order convergent,
at least when the frequency of large step ratios is very low ($N_1/N\approx 5\%$ in our tests).

\section{Concluding remarks}
\setcounter{equation}{0}

Consider some multi-step scheme
having the discrete kernels $\big\{B_{n-k}^{(n)}\big\}_{k=1}^n$ for parabolic equations,
\begin{align*}
\sum_{k=1}^nB^{(n)}_{n-k}\diff  u^k=\varepsilon\Delta u^n+f^n\quad\text{for $n\ge1$ and $B^{(n)}_{0}\neq0$.}
\end{align*}
We present a novel framework for the numerical analysis
by introducing a new class of DOC
kernels $\big\{\Theta_{n-k}^{(n)}\big\}_{k=1}^n$ defined via the orthogonal identity
\begin{align*}
\sum_{j=k}^{n}\Theta_{n-j}^{(n)}B^{(j)}_{j-k}\equiv\delta_{nk}
\quad\text{for $\forall\;1\leq k\le n$.}
\end{align*}
Taking the advantage of orthogonality, one has the following alternative form
\begin{align*}
\diff  u^n=\varepsilon\sum_{j=1}^n\Theta^{(n)}_{n-j}\Delta u^j+\sum_{j=1}^n\Theta^{(n)}_{n-j}f^j\quad\text{for $n\ge1$.}
\end{align*}
If the discrete kernels $B_{n-k}^{(n)}$ are positive semi-definite,
then the orthogonality implies the positive semi-definiteness of
DOC kernels $\Theta_{n-k}^{(n)}$.
So one has the following $L^2$ norm priori estimate
\begin{align*}
\mynormb{u^n}^2\le \mynormb{u^0}^2+2\sum_{k=1}^n\sum_{j=1}^k\myinnerb{u^k,\Theta^{(k)}_{k-j}f^j}\quad\text{for $n\ge1$.}
\end{align*}

For the adaptive BDF2 method applied to linear reaction-diffusion equations,
the above approach provides a concise stability and convergence theory,
which seems quite similar to that of the most robust BDF1 scheme.
Some applications will be reported in subsequent articles
for the numerical analysis of nonuniform BDF2 time-stepping scheme in simulating
the gradient flows \cite{LiaoJiZhang:2019,LiaoWangZhang:2019AC}, which always permit multiply time scales in approaching the steady state.
We expect that the novel theoretical framework will be useful to establish the optimal $L^2$ norm error estimate
for some other nonlocal time approximations having a discrete convolution form.

\section*{Acknowledgements}
The authors would like to thank Dr. Buyang Li for his valuable
discussions and fruitful suggestions.

\end{document}